\UseRawInputEncoding
\documentclass[11pt,twoside]{amsart}
\usepackage{amsmath}
\usepackage{amsthm}
\usepackage{amsfonts}
\usepackage{amssymb}
\usepackage[all]{xy}
\usepackage{alltt}
\usepackage{enumitem}
\usepackage{stmaryrd}
\usepackage{comment}
\usepackage{xspace}
\usepackage{enumitem}

\theoremstyle{plain}
\newtheorem{prop}{Proposition}[section]

\newtheorem{lem}[prop]{Lemma}

\newtheorem{thm}[prop]{Theorem}
\theoremstyle{definition}

\newtheorem{remark}[prop]{Remark}
\newtheorem{example}[prop]{Example}

\usepackage[margin=3.15cm]{geometry}

\usepackage[
        colorlinks=true,urlcolor=blue,linkcolor=blue,citecolor=blue
]{hyperref}

\newcommand{\rar}{\rightarrow}

\begin{document}
\title{Gaussian Binomial Coefficients in Group Theory, Field Theory, and Topology}

\author{Sunil K. Chebolu and Keir Lockridge}

\maketitle

\begin{abstract} In this article we offer group-theoretic, field-theoretic, and topological interpretations of the Gaussian binomial coefficients and their sum.
For a finite $p$-group $G$ of rank $n$, we show that the Gaussian binomial coefficient $\binom{n}{k}_p$ is the number of subgroups of $G$ that are minimally expressible as an intersection of $n - k$ maximal subgroups of $G$, and their sum is precisely the number of subgroups that are either $G$ or an intersection of maximal subgroups of $G$. We provide a field-theoretic interpretation of these quantities through the lens of Galois theory and a topological interpretation involving covering spaces.
\end{abstract}


\section{Introduction.}

The binomial coefficients $\binom{n}{k}$ arise naturally in a wide variety of combinatorial problems, but they appear perhaps most famously as the entries in Pascal's triangle, giving the coefficients in the binomial expansion \[(x + y)^n = \sum_{k = 0}^n \binom{n}{k} x^k y^{n-k}.\] Setting $x = y = 1$, one obtains an elegant formula for their sum: $2^n.$ Similarly, the sum of the multinomial coefficients in the expansion of $(x_1 + \cdots + x_m)^n$ is $m^n$.

The {\em Gaussian binomial coefficients} $\binom{n}{k}_q$ are defined by 
\[ \binom{n}{k}_q = \begin{cases} \frac{(q^n-1)(q^{n}-q)\cdots(q^{n}-q^{k-1})}{(q^k-1)(q^{k}-q)\cdots(q^k-q^{k-1})} \ \ \ \ \ \ &1 \le k \le n, \\
\ \ \ \ \ \ \ \ \ \ \ \ \ 1 & k = 0.
\end{cases}\] They too are the coefficients in an expansion of the expression $(x + y)^n$, but in a different setting: Assume $x$ and $y$ are noncommuting variables, and then impose a cost on commuting $x$ past $y$ with the relation $xy = qyx$. After expanding and collecting like terms using this relation, we have \[ (x + y)^n = \sum_{k = 0}^n \binom{n}{k}_q x^k y^{n-k}\] (see \cite{cohn}). Except for the case $q = 1$ (recovering the ordinary binomial coefficients), this equation does not provide a nice formula for the sum $S_{n, q} = \sum_{k=0}^n \binom{n}{k}_q$ because the polynomial evaluation map $f(x, y) \mapsto f(1, 1) \in \mathbb{Z}$ is not a ring homomorphism (it sends the zero polynomial $0 = xy - qyx$ to $1 - q$). The equation does show that the sum is congruent to $2^n$ modulo $q - 1$. As far as we know, there is no known simple formula for $S_{n, q}$.

There are, however, a number of interesting ways to interpret the Gaussian binomial coefficients and their sum. For instance, $\binom{n}{k}_q$
is the number of $k$-dimensional subspaces of $\mathbb{F}_q^n$, the $n$-dimensional vector space over the finite field $\mathbb{F}_q$. The number $S_{n, q}$ is therefore the total number of subspaces of $\mathbb{F}_q^n$. There are also combinatorial interpretations (involving lattice paths) and number-theoretic interpretations (involving integer partitions); see \cite{Stanley}.

In this article, we give group-theoretic, field-theoretic, and topological interpretations of the Gaussian binomial coefficients and their sum. Let $G$ be a finite group. A maximal subgroup of $G$ is a proper subgroup of $G$ that is not contained in any other proper subgroup of $G$. The intersection of all the maximal subgroups of $G$ is called the {\em Frattini subgroup} of $G$, denoted $\Phi(G)$. A natural question arises: \, {\em How many other subgroups of $G$ may be obtained by intersecting some number of maximal subgroups, and, more specifically, how many may be obtained by intersecting exactly $k$ maximal subgroups but not fewer?} \, Theorem \ref{mainthm} in Section 3 answers this question for finite $p$-groups in terms of the Gaussian binomial coefficients and their sum. This result uses the correspondence theorem in group theory between the lattice of subgroups of a group $G$ that contain a fixed normal subgroup $N$ and the lattice of subgroups in the quotient $G/N$. In Sections 4 and 5, we use two other well-known ``Galois correspondences" to translate this result into field theory and topology. The first is the correspondence between the subfield lattice of a Galois extension and the subgroup lattice of its Galois group; the second is the correspondence between the lattice of covering spaces of a sufficiently nice space and the subgroup lattice of its fundamental group.

\section{A Linear-algebraic interpretation.} In Section 3 we will see that a finite $p$-group whose Frattini subgroup is trivial is a finite vector space over $\mathbb{F}_p$, and the general case reduces to this situation. We will therefore first answer our question for $n$-dimensional $\mathbb{F}_p$-vector spaces. As mentioned in the introduction, it is well known that $\binom{n}{k}_q$ counts the number of subspaces of dimension $k$ in an $n$-dimensional vector space over the finite field $\mathbb{F}_q$. We will draw a conclusion after making just two further observations.

First, it is straightforward to check that counting subspaces of an $\mathbb{F}_p$-vector space is the same as counting subgroups: the two concepts coincide since $p$ is prime. (They do not coincide for vector spaces over finite fields in general.)

Second, it is well known that in a finite-dimensional vector space, every proper subspace is an intersection of codimension-1 (maximal) subspaces. Below, we state a version of this result that tracks the dimension of the intersections, and we include a proof for convenience. Maximal subspaces are called {\em linear hyperplanes}, and, according to our first observation, these are precisely the maximal subgroups in a finite-dimensional $\mathbb{F}_p$-vector space. By convention, the intersection of zero linear hyperplanes in a finite-dimensional vector space $V$ is $V$ itself.

\begin{lem}\label{allareints} Let $V$ be an $n$-dimensional vector space. A subspace of $V$ has dimension $k$ if and only if it is the intersection of $n-k$ (distinct) linear hyperplanes. In particular, every subspace of $V$ is an intersection of linear hyperplanes.
\end{lem}

\begin{proof}  
The vector space $V$ is the unique subspace of dimension $n$ and it is the intersection of zero linear hyperplanes, so assume $0 \leq k < n$ and let $\{v_1, \dots, v_k\}$ be a basis for $H$ (this basis is empty if $k = 0$). Complete this basis to a basis $B = \{v_1, \dots, v_k, v_{k+1}, \dots, v_n\}$ for $V$. Let $H_i$ denote the linear hyperplane spanned by all the vectors in $B$ except $v_i$. Then $H$ is the intersection of the $n-k$ linear hyperplanes $H_{k+1}, \dots, H_n$. To prove the converse, it suffices to prove that, for any linear hyperplane $H$ and subspace $K$ that is not contained in $H$, we have $\dim(H \cap K) = \dim K - 1$. Since $H \cap K$ is a proper subspace of $K$, we certainly have $\dim(H \cap K) \leq \dim K - 1$. But since there cannot be a linearly independent set with more than one vector contained in $V \setminus H$, we also have $\dim(H \cap K) \geq \dim K - 1$. This completes the proof.
\end{proof}

If $H$ is a subgroup of $G$ that is expressible as an intersection of maximal subgroups, then define the {\em intersection number} of $H$ to be the smallest positive integer $n$ such that $H$ is the intersection of $n$ maximal subgroups; by convention, the intersection number of $G$ is 0. We then have the following theorem (recall that subspaces and subgroups coincide for the vector space $V$ below).

\begin{thm} Let $V$ be an $n$-dimensional $\mathbb{F}_p$-vector space. \label{fpthm} \label{count}
\begin{enumerate}
\item The number of subspaces of $V$ that are either $V$ or an intersection of maximal subspaces is $S_{n, p}$.
\item The number of subspaces of $V$ with intersection number $n - k$ is $\binom{n}{k}_{p}.$
\item The number of subspaces of $V$ with dimension $k$ is $\binom{n}{k}_{p}.$
\end{enumerate}
\end{thm}

\noindent In the next section, we turn to the general case.

\section{A group-theoretic interpretation.} The {\em rank} of a group $G$ is the smallest positive integer $n$ such that $G$ can be generated by $n$ elements. We will soon see that, for a finite $p$-group $G$, the subgroups that are either $G$ or an intersection of maximal subgroups are precisely the subgroups containing the Frattini subgroup. This is not true in general: for example, the Frattini subgroup of the symmetric group $S_4$ is trivial, and none of its cyclic subgroups of order 4 are expressible as an intersection of maximal subgroups.

Let's now consider an example that illustrates the correspondence between subgroups of $p$-groups expressible as an intersection of maximal subgroups and proper subspaces of a finite vector space. 

\begin{example} \label{q8}
Consider $Q_8 = \{ \pm 1, \pm i, \pm j, \pm k\}$, the 8-element quaternion group with the relations $i^2 = j^2= k^2 = -1, ij = k, jk = i$ and $ki = j$. The quaternion group is a 2-group of rank 2 with exactly three maximal subgroups: $M_1 = \{ \pm 1, \pm i \}$, $M_2 = \{ \pm 1, \pm j \}$ and $M_3 = \{ \pm 1, \pm k \}$. Note that $M_1 \cap M_2 = M_2 \cap M_3 = M_3 \cap M_1 = M_1 \cap M_2 \cap M_3 =  \{ \pm 1 \}$, so $\Phi(Q_8) = \{\pm 1\}.$ Now, $Q_8/\{\pm 1\} \cong \mathbb{F}_2 \times \mathbb{F}_2$ is a 2-dimensional vector space over $\mathbb{F}_2$ with basis given by the images of $i$ and $j$ in the quotient. It has $S_{2,2}$ subspaces. The three linear hyperplanes $\mathbb{F}_2 \times \{0\}, \{0\} \times \mathbb{F}_2,$ and $\langle(1, 1)\rangle$ correspond to $M_1, M_2,$ and $M_3$, respectively; the trivial subspace corresponds to $\Phi(Q_8)$, and $\mathbb{F}_2 \times \mathbb{F}_2$ corresponds to $Q_8$ itself.

Thus, the number of subgroups of $Q_8$ that are expressible as an intersection of maximal subgroups or $Q_8$ is $S_{2,2} = 5$. Further, $\binom{2}{k}_2$ is the number of subgroups with intersection number $n-k$ because $\binom{2}{0}_2 = 1$, $\binom{2}{1}_2 = 3$, and $\binom{2}{2}_2 = 1$. We will see in Theorem \ref{mainthm} that this count depends only on the prime $p$ and the rank of the group, so the conclusion here would be the same for any 2-group of rank 2.
\end{example}

We now recall some basic facts about $p$-groups which may be found in \cite{DummitFoote}.
A finite $p$-group is a group whose order is a power of a prime $p$. Let $p$ be prime and let $G$ be a finite $p$-group of rank $n$. The Frattini subgroup $\Phi(G)$ is a normal subgroup of $G$, and $G/ \Phi(G)$ is an elementary abelian $p$-group of rank $n$; we therefore identify $G/\Phi(G)$ with $\mathbb{F}_{p}^n$. In fact, according to Burnside's basis theorem, $\{g_1, g_2, \ldots, g_n\}$ is a minimal generating set for $G$ if and only if $\{ \pi(g_1), \pi(g_2), \ldots, \pi(g_n)\}$ is an $\mathbb{F}_p$-basis for $\mathbb{F}_p^n$, where \[\pi \colon G \rar G/ \Phi(G)\cong \mathbb{F}_{p}^{n}\] is the quotient map.

\begin{lem}
Let $G$ be a finite $p$-group of rank $n$. \label{mainlem}
\begin{enumerate}
\item \label{one-to-one} The map $\pi$ induces a one-to-one correspondence between the subgroups of $G$ that contain $\Phi(G)$ and the subspaces of $\mathbb{F}_{p}^{n}$.
\item \label{codim} The map $\pi$ induces a one-to-one correspondence between the maximal subgroups of $G$ and the linear hyperplanes in $\mathbb{F}_{p}^{n}$. In particular, there are $(p^{n}-1)/(p-1)$ maximal subgroups in $G$.
\item \label{intsareall} Every proper subgroup of $G$ containing $\Phi(G)$ is an intersection of maximal subgroups.
\item \label{inumitem} For any subgroup $H$ with $\Phi(G) \leq H \leq G$, the intersection number of $H$ is $n-k$ if and only if $\pi(H)$ has dimension $k$.
\end{enumerate}
\end{lem}

\begin{proof} 
Item \ref{one-to-one} is a consequence of the correspondence theorem for groups and the fact that subgroups and subspaces coincide for $\mathbb{F}_p^n$.

Item \ref{codim} follows from the fact that $\pi$ must take maximal subgroups to maximal subgroups. Now, the number of maximal subgroups of $G$ is the number of linear hyperplanes in $\mathbb{F}_p^n$, which is $\binom{n}{n-1}_p = (p^{n}-1)/(p-1)$ by Lemma \ref{count}.

For item \ref{intsareall}, note first that every proper subspace of $\mathbb{F}_{p}^{n}$ is an intersection of linear hyperplanes by Lemma \ref{allareints}. The result now follows from item \ref{codim} and the elementary fact that $\pi^{-1}(\cap_{i} H_{i}) = \cap_{i} \pi^{-1}(H_{i})$, where the $H_{i}$ are linear hyperplanes in $\mathbb{F}_{p}^{n}$.

In light of Lemma \ref{allareints}, to prove item \ref{inumitem} it suffices to prove that $H$ and $\pi(H)$ have the same intersection number. This follows from item \ref{codim} and the elementary fact cited in the previous paragraph.
\end{proof}


By Lemma \ref{mainlem}, the subgroups of a finite $p$-group $G$ of rank $n$ that are an intersection of maximal subgroups or $G$ are in one-to-one correspondence with the subspaces of $\mathbb{F}_p^n$. Lemma \ref{count} therefore allows us to count them according to both their ranks modulo the Frattini subgroup and their intersection numbers. We now have the following group-theoretic interpretation of the Gaussian binomial coefficients and their sum.

\begin{thm} \label{group}
Let $G$ be a finite $p$-group of rank $n$. \label{mainthm}
\begin{enumerate}
\item The number of subgroups of $G$ that contain the Frattini subgroup is $S_{n, p}$.
\item The number of subgroups of $G$ with intersection number $n - k$ is $\binom{n}{k}_{p}.$

\item The number of subgroups of $G$ that contain the Frattini subgroup and whose dimension modulo the Frattini subgroup is $k$ is given by $\binom{n}{k}_{p}.$

\end{enumerate}
\end{thm}

\noindent Since the set of $k$-dimensional subspaces of a finite-dimensional vector space is in one-to-one correspondence with the set of $(n-k)$-dimensional subspaces (under vector space duality), $\binom{n}{k}_p = \binom{n}{n-k}_p.$ It follows that the number of subgroups of $G$ with intersection number $k$ is also $\binom{n}{k}_p.$ 

\begin{example}
Let $G$ be a finite $p$-group of rank $n$. The group $G$ has rank $1$ if and only if it is cyclic, if and only if there is a unique maximal subgroup. For groups of rank $n > 1$, when does every pair of maximal subgroups intersect in the Frattini subgroup? The total number of maximal subgroups of $G$ is $\binom{n}{1}_{p}$, and the number of subgroups that are an intersection of maximal subgroups is $S_{n,p} - 1$. We want \[S_{n,p} - 1 - \binom{n}{1}_{p} = 1 \iff \sum_{k=2}^n \binom{n}{k}_p = 1.\] This holds if and only if the $p$-group $G$ has rank $n= 2$.
\end{example}

\section{A field-theoretic interpretation.}
A famous theorem of Shafarevich (see \cite{shaf}) implies that every finite $p$-group is realizable as a Galois group over $\mathbb{Q}$; this inspired us to recast the group-theoretic interpretation above (Theorem \ref{mainthm}) into the field-theoretic interpretation given in this section. See \cite{DummitFoote} for background material.

Let $L/K$ be a Galois extension whose Galois group is a finite $p$-group with rank $n$. An intermediate field of the extension $L/K$ is a field $E$ such that $K \subseteq E \subseteq L$. The fundamental theorem of Galois theory gives an order-reversing correspondence between the lattice of subgroups of the Galois group $G = \mathrm{Gal}(L/K)$ and the lattice of intermediate fields $E$ of the extension $L/K$. Under this correspondence, Â if $H_1$ and $H_2$ are subgroups of $G$ that correspond respectively to subfields $F_1$ and $F_2$ of $L$, then $H_1 \cap H_2$ corresponds to $F_1F_2$ (the compositum of $F_1$ and $F_2$, which is the subfield generated by $F_1$ and $F_2$). Since our extension has degree $p^n$, a minimal intermediate field is an intermediate field $M$ such that $M/K$ has degree $p$. Let $K_{\phi}$ denote the intermediate subfield that corresponds to the Frattini subgroup $\Phi(G)$ of $G$. $K_{\phi}$ is the compositum of all the minimal intermediate fields and it is also the maximal elementary abelian $p$-extension of $K$ within $L$. If $F$ is an intermediate subfield that is a compositum of minimal intermediate fields, then define the {\em compositum number} of $F $ to the smallest number of minimal intermediate fields whose compositum is $F$. We then have the following theorem.

\begin{thm} \label{field} Let $L/K$ be a Galois extension such that $\text{Gal}(L/K)$ is a finite $p$-group of rank $n$.
\begin{enumerate}
\item The number of intermediate fields that are contained in $K_{\phi}$ is $S_{n,p}$. \label{gal1}
\item The number of intermediate fields that are contained in $K_{\phi}$ with compositum number $n-k$ is $\binom{n}{k}_{p}.$ \label{gal3}
\item The number of intermediate fields that are contained in $K_{\phi}$ such that $\text{Gal}(K_{\phi}/E)$ has rank $k$ is $ \binom{n}{k}_{p}$. \label{gal2}
\end{enumerate}
\end{thm}

\begin{proof} Under Galois correspondence, the intermediate fields that are contained in $K_{\phi}$ corresponds to subgroups of the Galois group that contain the Frattini subgroup, which, by part 1 of Theorem \ref{group}, is $S_{n,p}$. This establishes item \ref{gal1}.

Item \ref{gal3} follows from part 2 of Theorem \ref{group} and the Galois correspondence.

For item \ref{gal2}, let $E$ be an intermediate field and let $H$ be $\mathrm{Gal}(L/E)$. Since $\Phi(G)$ is a normal subgroup in $G = \mathrm{Gal}(L/K) $, it is also normal in the subgroup $\mathrm{Gal}(L/E)$. Thus, by the fundamental theorem of Galois theory, we know that $K_\phi /E$ is a Galois extension and $\mathrm{Gal}(K_\phi/E) = H/\Phi(G)$. The result now follows from part 3 of Theorem \ref{group}.
\end{proof}

\begin{remark} The hypothesis of the above theorem has a complete 
field-theoretic translation. The Galois group $G =\text{Gal}(L/K)$ is a finite $p$-group if and only if $L/K$ is a Galois extension whose degree is a power of $p$, and $G$ has rank $n$ if and only if $K_\phi/ K$ is an extension of degree $p^n$. To see this, note that $G$ has rank $n$ if and only if $G/\Phi(G)$ has dimension $n$ as an $\mathbb{F}_p$-vector space, or equivalently, $G/\Phi(G)$ is a group of order $p^n$. According to the Galois correspondence, $G/\Phi(G)$ is the Galois group of the Galois extension $K_\phi/K$. Thus, $G$ has rank $n$ if and only if $K_\phi/K$ has degree $p^n$.
\end{remark}

\begin{example} We briefly summarize the translation of Example \ref{q8} into the field-theoretic setting by giving an example of a field extension whose Galois group is $Q_8$. Consider the extension \[ L\,/\,K = \mathbb{Q}(i, \sqrt[8]{2})\,/\,\mathbb{Q}(\sqrt{2}\, i).\] It can be shown that $\mathrm{Gal}(L/K) = Q_8$; see \cite{DummitFoote}. It then follows from the Galois correspondence that there are only 3 minimal intermediate fields that correspond to the fixed fields of the three maximal subgroups $M_1, M_2,$ and $M_3$ of $Q_8$ (see Example \ref{q8}). These fields are $F_1 = \mathbb{Q}(i, \sqrt{2}), F_2 = \mathbb{Q}((i+1)\sqrt[4]{2}),$ and $F_3 = \mathbb{Q}((i-1)\sqrt[4]{2})$, the three quadratic extensions over $K$. One can check that \[ F_1F_2 = F_1F_3 = F_2F_3 = F_1F_2F_3 = K_\phi = \mathbb{Q}(i, \sqrt[4]{2})\] (the maximal elementary abelian $2$-extension). The number of subfields of $L$ that are either $K$ or a compositum of minimal intermediate fields is \[ |\{K, F_1, F_2, F_3, K_2 \}| = \binom{2}{0}_2 + \binom{2}{1}_2 + \binom{2}{2}_2 = 5.\]
\end{example}

\section{A topological interpretation.} In this section we present a topological analog of the previous section using covering spaces. A {\em covering space} of a space $X$ is a space $Y$ together with a special map $p\colon Y \longrightarrow X$ called a {\em covering map}: for every $x \in X$, there is an open neighborhood $U$ of $x$ in $X$ such that $p^{-1}(U)$ is a disjoint union of open sets that $p$ maps homeomorphically onto $U$. A key example is the map $E\colon \mathbb{R} \longrightarrow S^1$ defined by $E(t) = e^{2\pi it}$ that wraps the real number line around the unit circle. The symbol $Y$ alone is used to refer to the covering space, and the map is understood. See \cite{Hatcher, lee} for more background material. All spaces below will be base-pointed and all maps will be base-point preserving continuous functions.

Every finitely presented group occurs as the fundamental group of a compact, connected, smooth manifold \cite[Exercise 8.1.8.1]{stillwell}. In particular, this applies to all finite $p$-groups, so let $X$ be a compact, connected, smooth manifold such that its fundamental group $\pi_1(X)$ is a finite $p$-group of rank $n$. The conditions on $X$ guarantee the existence of a universal (simply connected) covering space for $X$, and we may apply the necessary results from covering space theory.

There is a partial order on the set of all isomorphism classes of path-connected covering spaces of $X$, defined as follows. Suppose $C_1$ and $C_2$ are path-connected covering spaces with covering maps $p_1 \colon C_1 \rightarrow X$ and $p_2 \colon C_2 \rightarrow X$. We say that $C_1 \succ C_2$ ($C_1$ {\em covers} $C_2$) if there is a map $f \colon C_1 \rightarrow C_2$ such that $p_1 = p_2 \circ f$. The join $C_1 \vee C_2$ is a path-connected covering space $C$ with $C \succ C_1$ and $C \succ C_2$ such that if $C'$ is any other covering space with $C' \succ C_1$ and $C' \succ C_2$, then $C' \succ C$.  
 The fundamental theorem of covering spaces gives an order-reversing bijection between the lattice of subgroups of $\pi_1(X)$ and the lattice of isomorphism classes of path-connected covering spaces of $X$.  Under this bijection, the intersection of subgroups corresponds to the join of covering spaces, the trivial subgroup to the universal cover, the maximal subgroups to the minimal covering spaces, and normal subgroups to regular covers (see \cite{Hatcher}). 
 
 We will call the isomorphism class of the cover that corresponds to the Frattini subgroup of $\pi_1(X)$ the {\em Frattini cover.}
 Denote by $\pi_{X'}$ the subgroup of $\pi_1(X)$ that corresponds to a path-connected covering space $X'$ of $X$. Define the {\em join number} of a covering space $C$ that is a join of minimal coverings to be the smallest number of minimal coverings whose join is $C$. We now have the following analog of Theorem \ref{group}.

 \begin{thm}
 Let $X$ be a compact, connected, smooth manifold such that its fundamental group $\pi_1(X)$ is a finite $p$-group of rank $n$.
 \begin{enumerate}
\item The number of isomorphism classes of path-connected covering spaces of $X$ that are covered by the Frattini cover is $S_{n,p}$.
 \item The number of isomorphism classes of path-connected covering spaces of $X$ that are covered by the Frattini cover with join number $n-k$ is $\binom{n}{k}_{p}$.
\item The number of isomorphism classes of path-connected covering spaces $X'$ of $X$ that are covered by the Frattini cover and such that $\pi_{X'}/\Phi(\pi_1(X))$ has rank $k$ is $ \binom{n}{k}_{p}$.
\end{enumerate}
\end{thm}

\begin{example}
The projective plane $ \mathbb{RP}^2 $ is obtained as a quotient space of the sphere $S^2$ by identifying pairs of antipodal points; let $q \colon S^2 \rightarrow \mathbb{RP}^2$ denote the corresponding quotient map. It has fundamental group $\pi_1( \mathbb{RP}^2 ) = C_2$, a cyclic group of order two generated by $\sigma$. The space $X = \mathbb{RP}^2 \times \mathbb{RP}^2$ is a compact 4-dimensional connected smooth manifold with $\pi(X) = \pi_1( \mathbb{RP}^2 ) \times \pi_1( \mathbb{RP}^2 ) = C_2 \times C_2$, a 2-group of rank 2. The universal covering space of $X$ is given by the natural quotient map $ q \times q \colon \tilde{X} = S^2 \times S^2 \rightarrow \mathbb{RP}^2 \times \mathbb{RP}^2$. Since the Frattini subgroup is trivial, the Frattini covering space is the universal covering space. The 5 subgroups of $C_2 \times C_2$ act on $\tilde{X}$, and all path-connected covering spaces of $X$ arise as quotients of these group actions. The universal cover and the trivial cover correspond to the trivial group and the whole group, respectively. The three nontrivial (minimal) covers corresponding to the three nontrivial (maximal) subgroups of $C_2 \times C_2$ are:
 \[\begin{aligned}
 \tilde{X}/\langle 1, \sigma \rangle &\cong  S^2 \times \mathbb{RP}^2 \overset{1 \times q }\longrightarrow  \mathbb{RP}^2 \times \mathbb{RP}^2 \\
 \tilde{X}/\langle \sigma, 1 \rangle &\cong  \mathbb{RP}^2 \times S^2 \overset{q \times 1 }\longrightarrow \mathbb{RP}^2 \times \mathbb{RP}^2 \\
 \tilde{X}/\langle \sigma, \sigma \rangle &\cong (S^2\times S^2)/\{(x, y) \sim (-x, -y) \} \overset{\rho}\longrightarrow \mathbb{RP}^2 \times \mathbb{RP}^2 
 \end{aligned}\]
 where $\rho (\pm(u, v)) = ([\pm u], [\pm v])$.  The join of any two of these three covers is the universal cover.  Indeed, the number of covering spaces of $X$ that are either trivial or a join of minimal covering spaces is 
 \[ S_{2,2} = \binom{2}{0}_2 + \binom{2}{1}_2 + \binom{2}{2}_2 = 1+3+1= 5.\]
 
\end{example}

\noindent
\textbf{Acknowledgments.}
We thank Jon Carlson for providing helpful feedback on an early version of this article. We also thank  Andy Schultz and Craig Westerland for their suggestions.

\bibliographystyle{alpha}

\end{document}